\documentclass[journal]{IEEEtran}
\usepackage{graphicx}
\usepackage{amsmath}
\usepackage[justification=centering]{caption}
\usepackage{amssymb}

\begin{document}

\title{A Subspace Method for Array Covariance Matrix Estimation}

\author{Mostafa~Rahmani and George~K.~Atia,~\IEEEmembership{Member,~IEEE,} % <-this % stops a space
\thanks{This work was supported in part by NSF Grant CCF-1320547.

The authors are with the Department of Electrical Engineering and Computer Science, University of Central Florida,Orlando, FL 32816 USA (e-mail: mostafa@knights.ucf.edu, george.atia@ucf.edu).}% <-this % stops a spa
}

\markboth{}%
%\markboth{Journal of \LaTeX\ Class Files,~Vol.~11, No.~4, %December~2012}%
{Shell \MakeLowercase{\textit{et al.}}: Bare Demo of IEEEtran.cls for Journals}
% make the title area
\maketitle

% As a general rule, do not put math, special symbols or citations
% in the abstract or keywords.
\begin{abstract}
This paper introduces a subspace method for the estimation of an array covariance matrix. It is shown that when the received signals are uncorrelated, the true array covariance matrices lie in a specific subspace whose dimension is typically much smaller than the dimension of the full space. Based on this idea, a subspace based covariance matrix estimator is proposed. The estimator is obtained as a solution to a semi-definite convex optimization problem. While the optimization problem has no closed-form solution, a nearly optimal closed-form solution is proposed making it easy to implement. In comparison to the conventional approaches, the proposed method yields higher estimation accuracy because it eliminates the estimation error which does not lie in the subspace of the true covariance matrices. The numerical examples indicate that the proposed covariance matrix estimator can significantly improve the estimation quality of the covariance matrix. 
\end{abstract}

% Note that keywords are not normally used for peerreview papers.

\IEEEpeerreviewmaketitle

\section{Introduction}
% and "HIS" in caps to complete the first word.
\IEEEPARstart{T}{he} 
 estimation of covariance matrices is a crucial component of many signal processing algorithms [1-4]. In many applications, there is a limited number of snapshots and the sample covariance matrix cannot yield the desired estimation accuracy. This covariance matrix estimation error significantly degrades the performance of such algorithms. 
In some applications, the true covariance matrix has a specific structure. For example, the array covariance matrix of a linear array with equally spaced antenna elements is a Toeplitz matrix when the sources are uncorrelated [5, 6]. Moreover, in some applications [4, 8], the structure of the problem suggests that the underlying true covariance matrix is the Kronecker product of two valid covariance matrices [4, 7].\par
This side information can be leveraged in covariance matrix estimation to improve the estimation quality. For instance, in [5] a weighted least square estimator for covariance matrices with Toeplitz structures was proposed and it was shown that the resulting covariance matrix can enhance the performance of angle estimation algorithms, such as MUltiple SIgnals Classification (MUSIC) [13]. In [8], covariance matrices with Kronecker structure are investigated and a maximum likelihood based algorithm is introduced.\par
In addition, the structure of covariance matrices has been exploited in various DOA estimation algorithms, such as the linear structure in [9], and the diagonal structure for the covariance matrix of uncorrelated signals in [10]. Recently some research works have focused on the application of sparse signal processing in DOA estimation based on the sparse representation of the array covariance matrix. For example, [11] proposes the idea that the eigenvectors of the array covariance matrix have a sparse representation over a dictionary constructed from the steering vectors. In [12, 14], it is shown that when the received signals are uncorrelated, the array covariance matrix has a sparse representation over a dictionary constructed using the atoms, i.e. the correlation vectors. A similar idea is proposed in [15], with the difference that the proposed method does not require choosing a hyper-parameter.\par
In this paper, we focus on the estimation of  array covariance matrices with linear structure. First, we show that when the sources are uncorrelated, the array covariance matrix has a linear structure implying that all possible array covariance matrices can be described by a specific subspace. Based on this idea, a subspace-based covariance matrix estimator is proposed as a solution to a semi-definite convex optimization problem. Furthermore, we propose a nearly optimal closed-form solution for the proposed covariance matrix estimator. Our results show that the proposed method can noticeably improve the covariance matrix estimation quality. Moreover, the closed-form solution is shown to closely approach the optimal performance.  

\section{ARRAY SIGNAL MODEL}

The system model under consideration is a narrowband array system with $N$ antennas. All the signals are assumed to be narrowband with the same center frequency and impinge on the array from the far field. The baseband array output can be expressed as
\begin{eqnarray}
\textbf{x}(t)= \sum\limits_{i=1}^p z_i(t)\textbf {v}(\theta_i,\phi_i)+\textbf {n}(t)
\end{eqnarray}
where $\textbf{x}(t)$ is the $N\times1$ array output vector, $p$ is the number of the received signals, $z_i(t)$ is the $i^{th}$ signal, $(\theta_i,\phi_i)$ is the elevation and azimuth arrival angle of the $i^{th}$ signal, $\textbf {v}(\theta_i,\phi_i)$ is the baseband array response to $i^{th}$ signal and \textbf {n}(t)  is the noise vector. The baseband array response, $\textbf {v}(\theta_i,\phi_i)$, is called the ``steering vector'' [13].\par
If the received signals are uncorrelated, the covariance matrix can be written as
\begin{eqnarray}
\textbf {R}= \sum\limits_{i=1}^p \sigma _i^2 \textbf {v}(\theta_i,\phi_i)  \textbf {v}^H(\theta_i,\phi_i)+\sigma _n^2\textbf{I}
\end{eqnarray}
where $\sigma _i^2$  represents the power of the $i^{th}$ signal, $\sigma _n^2$ is the noise variance and $\textbf{I}$ is the identity matrix.
\section{PROPOSED ALGORITHM}
We define the ``correlation vector'' which belongs to direction $(\theta,\phi)$   as follows
\begin{eqnarray}
\textbf{c}(\theta,\phi)=vec\displaystyle(\textbf{v}(\theta,\phi) \textbf{v}^H(\theta,\phi)) 
\end{eqnarray}
where $vec(\bullet)$ is a linear transformation that converts its matrix argument to a vector by stacking the columns of the matrix on top of one another. Consequently, the covariance matrix can be rewritten as
\begin{eqnarray}
vec( \textbf{R}-\sigma_n^2 \textbf{I})= \sum\limits_{i=1}^p \sigma_i^2 \textbf{c}(\theta_i,\phi_i)
\end{eqnarray}
Therefore, $vec( \textbf{R}-\sigma_n^2 \textbf{I})$ is a linear combination of the correlation vectors of the received signals.\par
According to (4), $vec( \textbf{R}-\sigma_n^2 \textbf{I})$ lies in the subspace of the correlation vectors. Hence, if we build the subspace spanned by all possible correlation vectors$\{\textbf{c}(\theta ,\phi)|  0\le \theta \le \pi ,  0\le \phi \le 2\pi\}$ , then $vec( \textbf{R}-\sigma_n^2 \textbf{I})$ completely lies in this subspace. For many array structures, the matrix $ (\textbf{v}(\theta,\phi) \textbf{v}^H(\theta,\phi))$ inherits some symmetry properties.
Accordingly, the correlation vectors cannot span an $ N^2 $ dimensional space. For example, when the incoming signals are uncorrelated, the covariance matrix of a uniform linear array is a Toeplitz matrix [5]. It is easy to show that all the $N\times N$ Toeplitz matrices can be described by a $2N-1$  dimensional space. \par
The subspace of the correlation vectors $\{\textbf{c}(\theta,\phi)|0\le \theta \le \pi , 0 \le \phi \le 2\pi \}$  can be obtained by constructing a positive definite matrix
\begin{eqnarray}
\textbf{S}=\int\limits_{0}^{2\pi} \int\limits_ {0}^\pi \textbf{c} (\theta,\phi) \textbf{c}^H (\theta,\phi)d\theta d\phi
\end{eqnarray}
where (5) is an element-wise integral. Based on (5), the subspace dimension of the correlation vectors  $\{\textbf{c}(\theta,\phi)|0\le\theta\le\pi , 0 \le \phi \le 2\pi \}$    is equal to the number of non-zero eigenvalues of the matrix $\textbf{S}$. Consequently, the subspace of the correlation vectors can be constructed using the eigenvectors which correspond to the non-zero eigenvalues.\\
Fig. 1 shows the eigenvalues of $\textbf{S}$ for a square planar array with 16 elements (the horizontal and vertical space between the elements is half a wavelength). One can observe that the number of non-zero eigenvalues is equal to 49. Therefore, for this array, the subspace of the correlation vectors can be constructed from the 49 eigenvectors corresponding to the non-zero eigenvalues. Note that for a 16-element linear array, we observe 31 non-zeros eigenvalues because the covariance matrix is a Toeplitz matrix[5]. For some array structures such as circular array, we may not observe zero eigenvalues but our investigation has shown that the subspace of the correlation vectors can be effectively approximated using the dominant eigenvectors (the eigenvectors corresponding to the dominant eigenvalues).
\begin{figure}[t!]
    \includegraphics[width=0.5\textwidth]{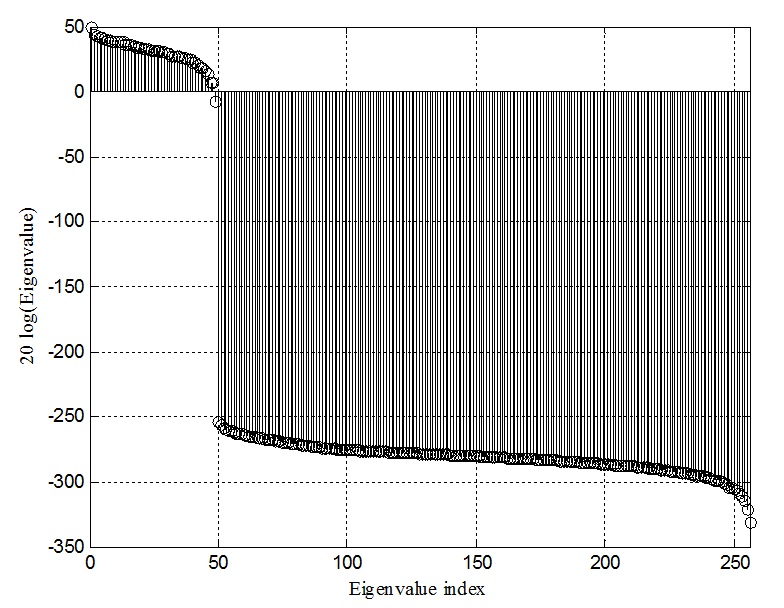}
    \caption{Eigenvalues of the matrix S}
    \centering
\end{figure}
Therefore, if we construct the matrix $\textbf{Q}$ whose columns form a basis for the correlation vectors subspace, we can rewrite the covariance matrix as 
\begin{eqnarray}
vec(\textbf{R}-\sigma _n^2\textbf{I})=\textbf{Q}\textbf{a} .
\end{eqnarray}
Hence, we can choose the columns of $\textbf{Q}$  as the eigenvectors corresponding to the non-zero eigenvalues (or the dominant eigenvectors). By imposing the linear structure constraint (6) to the covariance matrix estimation problem, we can significantly improve the estimation quality. Some works have studied covariance matrices with linear structures. For example, a weighted least-square estimator was proposed in [5] based on the linear structure for Toeplitz covariance matrices. However, the Toeplitz structure is restricted to linear arrays and the resulting matrix is not guaranteed to be positive definite. 
\subsection{Subspace Based Covariance Matrix Estimation}
Based on (4) and (6), the estimated covariance matrix should lie in the subspace spanned by the columns of $\textbf{Q}$ . We are going to estimate $\textbf{R}_s$  which is defined as
\begin{eqnarray}
\textbf{R}_s=\textbf{R}-\sigma_n^2\textbf{I}
\end{eqnarray}
Based on the previous discussion, we propose the following optimization problem
\begin{eqnarray}
\begin{aligned}
& \underset{\textbf{R}_s}{\text{min}}
& & \parallel(\hat{\textbf{R}}-\sigma_n^2\textbf{I})-\textbf{R}_s\parallel_2^2\\
& \text{subject to}
& & (\textbf{I}-\textbf{Q}(\textbf{Q}^H\textbf{Q})^{-1}\textbf{Q}^H)vec(\textbf{R}_s)=0\\
& &  & \textbf{R}_s\succeq0\\
\end{aligned}
\end{eqnarray}
where
\begin{eqnarray}
\hat{\textbf{R}}=\frac{1}{M}\sum\limits_{m=1}^M\textbf{x}(m)\textbf{x}^H(m)
\end{eqnarray}
is the sample covariance matrix and $M$  is the number of time samples. The matrix $(\textbf{I}-\textbf{Q}(\textbf{Q}^H\textbf{Q})^{-1}\textbf{Q}^H)$ is the projection matrix on the subspace that is orthogonal to the subspace of the correlations vectors. As such, the first constraint in (8) ensures that the resulting matrix lies in the correlations vectors subspace. The second constraint guarantees that the resulting matrix is positive definite.
Note that, (8) is a convex optimization problem and can be solved using standard tools from convex optimization. The proposed method imposes the linear structure using the subspace constraint. If the covariance matrix is Toeplitz, the subspace constraint enforces the resulting matrix to be Toeplitz. However, the proposed algorithm is not limited to Toeplitz structures and can be used for any linear structure.\par
The sample covariance matrix (9) can be expressed as
\begin{eqnarray}
\hat{\textbf{R}}=\sum\limits_{i=1}^P\sigma_i^2\textbf{v}(\theta_i)\textbf{v}^H(\theta_i)+\boldsymbol{\Delta}+\sigma_n^2\textbf{I}.
\end{eqnarray} 
The second term on the right hand side of (10), $\boldsymbol{\Delta} $ , is the unwanted part (estimation error) which tends to zero if we have an infinite number of snapshots. The estimation error has some random behavior and can lie anywhere in the entire space. Since the first constraint in (8) enforces the estimated matrix to lie in the correlation vectors' subspace, it is expected to eliminate the component of estimation error which is not in this subspace. The dimension of the correlation vectors subspace is typically smaller than the entire space dimension $N^2$. For example, for a 30-element uniform linear array, the dimension of the correlation vectors subspace is equal to 59; while the entire space dimension is 900. Thus, it is conceivable that the proposed method could yield a much better estimation performance in comparison to the sample covariance matrix (9). 
\subsection{Near Optimal Closed-form Solution}
The proposed optimization problem (8) is an $N^2$   dimensional optimization problem. Therefore, it may be hard to solve for large arrays. In this section, we derive a closed form near optimal solution which makes our method easy for practical implementation. \par
According to (4), the covariance matrix should be in the correlation vectors subspace. We define $\textbf{R} _ \perp$   and $\textbf{R} _ \parallel$  as follows
\begin{eqnarray}
vec(\textbf{R}_\parallel)=& \textbf{Q}(\textbf{Q}^H\textbf{Q})^{-1}\textbf{Q}^Hvec(\hat{\textbf{R}}-\sigma_n^2\textbf{I})\\
vec(\textbf{R}_\perp)=&(\textbf{I}-\textbf{Q}(\textbf{Q}^H\textbf{Q})^{-1}\textbf{Q}^H)vec(\hat{\textbf{R}}-\sigma_n^2\textbf{I}). 
\end{eqnarray}
Thus, $\textbf{R}_\perp$  is orthogonal to the correlation vectors subspace and $\textbf{R}_\parallel$  contains the desired part. Therefore, we rewrite (8) as
\begin{eqnarray}
\begin{aligned}
& \underset{\textbf{R}_s}{\text{min}}
& & \parallel\textbf{R}_s-\textbf{R}_\parallel\parallel_2^2\\
& \text{subject to}
& & (\textbf{I}-\textbf{Q}(\textbf{Q}^H\textbf{Q})^{-1}\textbf{Q}^H)vec(\textbf{R}_s)=0\\
& &  & \textbf{R}_s\succeq0\\
\end{aligned}
\end{eqnarray}
In the proposed estimator (8), we placed the first constraint to suppress the estimation error which does not lie in the correlation vectors subspace. In (11), we project the sample covariance matrix to the correlation vectors subspace. Thus, we have eliminated the estimation error which does not lie in the correlation vectors subspace. Accordingly, we simplify (13) as follows
\begin{eqnarray}
\begin{aligned}
& \underset{\textbf{R}_s}{\text{min}}
& & \parallel\textbf{R}_s-\textbf{R}_\parallel\parallel_2^2\\
& \text{subject to}
& & \textbf{R}_s\succeq0\\
\end{aligned}
\end{eqnarray}
which has a simple closed form solution
\begin{eqnarray}
\hat{\textbf{R}}_s=(\sum\limits_{i=1}^q\lambda_i\boldsymbol{\beta}_i\boldsymbol{\beta}_i^H)
\end{eqnarray}
where $q$   is the number of positive eigenvalues of $\textbf{R}_\parallel$ , $ \{\lambda_i\}_{i=1}^q$  are the positive eigenvalues and $ \{\boldsymbol{\beta}_i\}_{i=1}^q$   are their corresponding eigenvectors. Actually, we break the primary optimization problem (8) into two optimization problems. First, we find a matrix in the correlation vectors subspace which is most close to the sample covariance matrix and the resulting matrix is $\textbf{R}_\parallel$. In the second step, we find the closest positive semi-definite matrix to $\textbf{R}_\parallel$ and the resulting matrix is given in (15).

\section{Simulation Results}
In this section, we provide some simulation results to illustrate the performance of the proposed approach. The examples provided include DoA estimation and subspace estimation, which underscores the flexibility of the proposed covariance matrix estimation approach for a broad range of applications. All the curves are based on the average of 500 independent runs.
\subsection{Simulation I (DOA estimation, probability of resolution)}
Assume a uniform linear array with $N=10$ omnidirectional sensors spaced half a wavelength apart. For this array the correlation vectors subspace is a 19 dimensional space since the covariance matrix is Toeplitz. The additive noise is modeled as a complex Gaussian zero-mean spatially and temporally white process with identical variances in each array sensor. In this experiment, we compare the performance of MUSIC when used with the sample covariance matrix and with the proposed covariance matrix estimation method. We also compare its performance with the sparse covariance matrix representation method [14, 12] and the SParse Iterative Covariance-based Estimation approach (SPICE) [15]. We consider two uncorrelated sources located at $45^\circ$ and $50^\circ$ ($90^\circ$ is the direction orthogonal to the array line) and both sources are transmitted with the same power. Fig. 2 shows the probability of resolution (the probability that the algorithm can distinguish these two sources) versus the number of snapshots for one fixed sensor with $\mathrm{SNR}=0$ dB. It is clear that using the proposed method leads to significant improvement in performance in comparison to using the sample covariance matrix. SPICE [15] is an iterative algorithm, which is based on the sparse representation of the array covariance matrix and requires one matrix inversion in each iteration. One can see that this algorithm fails when we use 20 iterations, however, performs well with 1000 iterations. Nevertheless, for practical purposes it is generally computationally prohibitive to perform 1000 matrix inversion operations. In addition, one can observe that the proposed near optimal solution to (8) yields a close performance to the optimal solution. Fig. 3 displays the probability of resolution against $\mathrm{SNR}$ for a fixed training size $M=500$ snapshots. Roughly, the MUSIC algorithm based on the proposed method is 7 dB better than the MUSIC algorithm based on the sample covariance matrix. In summary, the proposed method yields notable and promising performance even with a small number of snapshots and at low SNR regimes. Furthermore, it is easily implementable using the proposed closed-form solution, which consists of a matrix multiplication and eigen-decomposition.
\begin{figure}[t!]
    \includegraphics[width=0.5\textwidth]{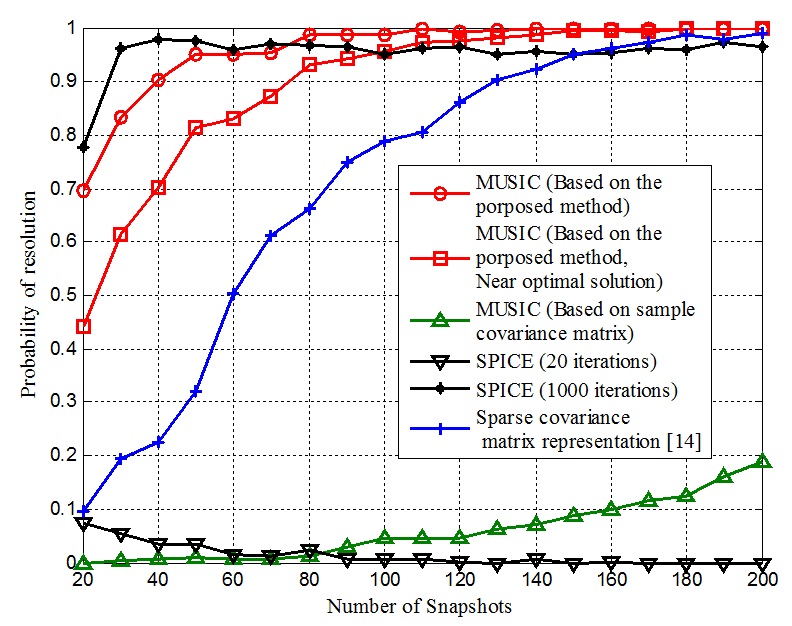}
    \caption{Probability of resolution as a function of the number of snapshots}
    \centering
\end{figure}

\begin{figure}[t!]
    \includegraphics[width=0.5\textwidth]{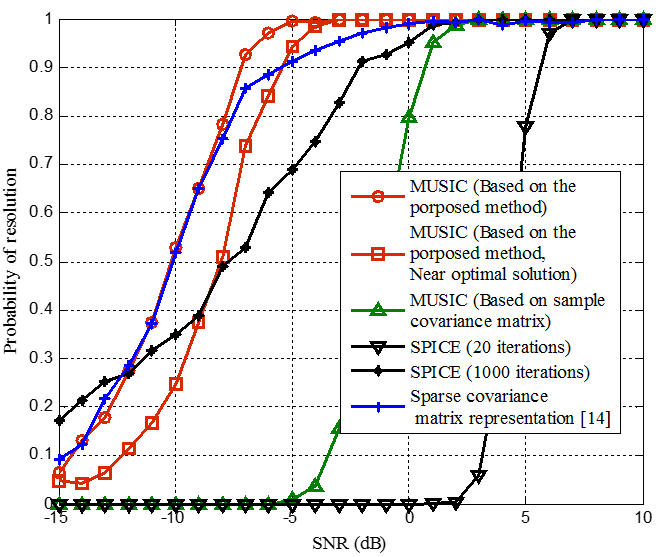}
    \centering
    \caption{Probability of resolution versus $\mathrm{SNR}$}
   
\end{figure}

\subsection{Simulation II (Signals subspace estimation)}
The estimation of the subspace of the received signals is an important task in many signal processing algorithms. For example, in the eigen-space based beamforming algorithm [3], the subspace of the received signals is used to make the beamformer robust against the steering vector mismatch. In the MUSIC algorithm, the subspace of the received signals is used to obtain the noise subspace [13]. The subspace of the received signals is usually estimated using the dominant eigenvectors of the estimated covariance matrix. In this simulation, we consider three uncorrelated sources located at $85^\circ$, $90^\circ$ and $95^\circ$ and the sources are received with same signal to noise ratio. To investigate the accuracy of the subspace estimation, we define the distance between two subspaces as follows [16]:\\
Given two matrices  $\hat{\textbf{U}}$,$\hat{\textbf{V}}\in\mathbb{R}^{N\times K}$, the distance between the subspaces spanned by the columns of $\hat{\textbf{U}}$ and $\hat{\textbf{V}}$ is defined as
\begin{eqnarray}
dist(\hat{\textbf{U}},\hat{\textbf{V}})=\parallel\textbf{U}^H_\perp \textbf{V} \parallel_2 = \parallel\textbf{V}^H_\perp \textbf{U} \parallel_2
\end{eqnarray}
where $\textbf{U}$ and $\textbf{V}$ are orthonormal bases of the spaces $span(\hat{\textbf{U}})$ and  $span(\hat{\textbf{V}})$, respectively. Similarly, $\textbf{U}^H_\perp\in\mathbb{R}^{N\times (N-K)} $ is an orthonormal basis for the subspace which is orthogonal to $span(\hat{\textbf{U}})$ and $\textbf{V}^H_\perp\in\mathbb{R}^{N\times (N-K)} $ is an orthonormal basis for the subspace which is orthogonal to $span(\hat{\textbf{V}})$. In addition, $\parallel\textbf{X}\parallel_2$ denotes the spectral norm of matrix $\textbf{X}$. Fig. 4 displays the distance between the true signals subspace and the estimated one as a function of the number of snapshots for $\mathrm{SNR}=-6$ dB. We construct the signal subspace using the first three eigenvectors. One can observe that the proposed method exhibits a better rate of convergence. In addition, the performance of the closed-form solution closely approaches the optimal solution. 
\begin{figure}[t!]
    \includegraphics[width=0.5\textwidth]{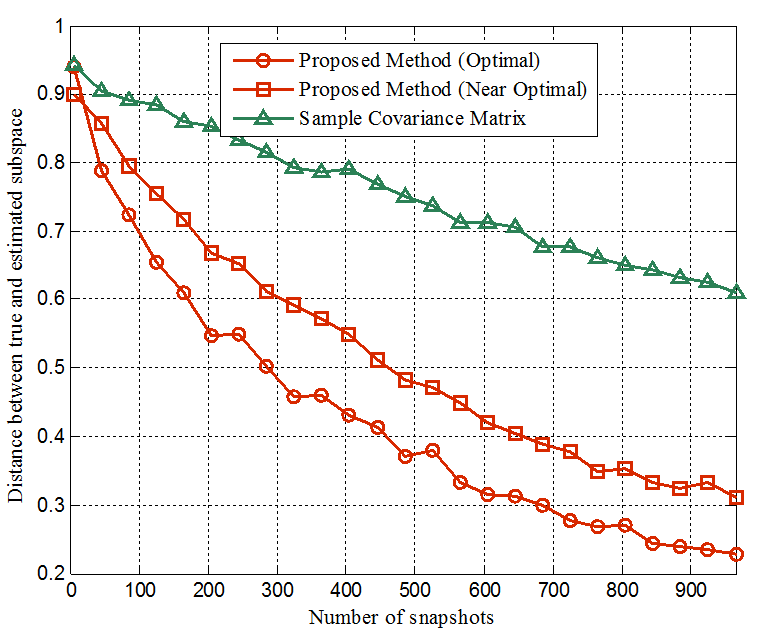}
    \centering
    \caption{Distance between the true and the estimated signals subspace}
   
\end{figure}
\section{Conclusion}
In this paper, a subspace method for array covariance matrix estimation was proposed. We have shown that when the received signals are uncorrelated, the covariance matrix lies in the subspace of the correlation vectors. Based on this idea, we posed the estimation problem as a convex optimization problem and enforced a subspace constraint. In addition, a near optimal closed-form solution for the proposed optimization problem was derived. A number of numerical examples demonstrated the notable performance of the proposed approach and its applicability to a wide range of signal processing problems, including but not limited to, DoA estimation and subspace estimation. In contrast to some of the existing approaches, which suffer from drastic performance degradation with limited data and at low SNR regimes, the proposed method showed very graceful degradation in such settings. 
\newpage
% references section

% can use a bibliography generated by BibTeX as a .bbl file
% BibTeX documentation can be easily obtained at:
% http://www.ctan.org/tex-archive/biblio/bibtex/contrib/doc/
% The IEEEtran BibTeX style support page is at:
% http://www.michaelshell.org/tex/ieeetran/bibtex/
%\bibliographystyle{IEEEtran}
% argument is your BibTeX string definitions and bibliography database(s)
%\bibliography{IEEEabrv,../bib/paper}
%
% <OR> manually copy in the resultant .bbl file
% set second argument of \begin to the number of references
% (used to reserve space for the reference number labels box)

\end{document}